\documentclass[a4paper,12pt]{article}
\usepackage{comment}
\usepackage{cite}
\usepackage{amsmath}
\usepackage{amssymb}
\usepackage{amsfonts}
\usepackage[T1]{fontenc}
\usepackage[utf8]{inputenc}
\usepackage{graphicx}
\usepackage{fancyhdr}
\usepackage{float}
\usepackage{xcolor}
\usepackage{authblk}
\usepackage{mathrsfs}
\usepackage{empheq}
\usepackage{url}

\usepackage{geometry}
\begin{document}

\title{From Fourier series to infinite product representations of $\pi$ and infinite-series forms for its positive powers}

\author[$\dagger$]{Jean-Christophe {\sc Pain}\footnote{jean-christophe.pain@cea.fr}\\
\small
CEA, DAM, DIF, F-91297 Arpajon, France\\
Universit\'e Paris-Saclay, CEA, Laboratoire Mati\`ere en Conditions Extr\^emes,\\ 
91680 Bruy\`eres-le-Ch\^atel, France
}

\maketitle

\begin{abstract}
In this article, we derive, using Fourier series and multiple derivative of the function $\pi/\sin(\pi x)$, series representations for positive powers of $\pi$. We also show that the Euler-Wallis product can be easily obtained from the same formalism and deduce infinite products representations of $\pi$.
\end{abstract}

\section{Introduction}

Let us set $f(x)=\cos(\alpha x)$. The Fourier transform $\mathscr{F}$ of $f$ satisfies $\mathscr{F}(x)=f(x)$ in $[-\pi,\pi]$, $\mathscr{F}(x+2\pi)=\mathscr{F}(x)$, $\mathscr{F}$ is an odd function, $\mathscr{C}^{\infty}$ and piecewise continuous. One has \cite{Tolstov1962,Kaplan1992}:
\begin{equation}
\mathscr{F}(x)=\frac{a_0}{2}+\sum_{n=1}^{\infty}a_n~\cos(nx),
\end{equation}
where
\begin{equation}
a_n=\frac{2}{\pi}\int_0^{\pi}\cos(\alpha t)\cos(n t)dt=\frac{(-1)^n}{\pi}\sin(\alpha\pi)\left(\frac{1}{\alpha+n}+\frac{1}{\alpha-n}\right).
\end{equation}
Setting $x=0$ and $x=\pi$, one gets respectively
\begin{empheq}[box=\fbox]{align}
\frac{\pi}{\sin(\alpha\pi)}=\sum_{n=-\infty}^{\infty}\frac{(-1)^n}{\alpha+n}
\end{empheq}
and
\begin{empheq}[box=\fbox]{align}\label{cota}
\pi~\mathrm{cotan}(\alpha\pi)=\sum_{n=-\infty}^{\infty}\frac{1}{\alpha+n}.
\end{empheq}

\section{Infinite series representation for $\pi^k$ with $k\ge 0$ involving multinomial coefficients}

Starting from
\begin{equation}
\frac{\pi}{\sin(\pi x)}=\sum_{n=-\infty}^{\infty}\frac{(-1)^n}{x+n},
\end{equation}
with $x=1/4$, one obtains
\begin{equation}
\pi=2\sqrt{2}~\left[1+\sum_{n=-\infty}^{\infty}\frac{(-1)^n}{(1+4n)}\right].
\end{equation}
Using the F\`aa di Bruno formula for the multiple derivative of a composite function (and in particular here of the inverse of a function):
\begin{equation}
\left(\frac{1}{f(x)}\right)^{(k)}=\frac{k!}{\left[f(x)\right]^{k+1}}\sum_{\substack{p_0+p_1+p_2+\cdots +p_k=k\\p_1+2p_2+\cdots+kp_k=k}}\mathscr{C}(p_1,p_2,\cdots,p_k)\prod_{i=0}^k\left[f^{(i)}(x)\right]^{p_i},
\end{equation}
with
\begin{equation}
\mathscr{C}(p_1,p_2,\cdots,p_k)=\frac{(-1)^{k-p_0}(k-p_0)!}{\prod_{i=1}^k(i!)^{p_i}p_i!}.
\end{equation}
$\mathscr{C}(p_1,p_2,\cdots,p_k)$ can be expressed in terms of the multinomial coefficient
\begin{equation}
\mathscr{M}(p_1,p_2,\cdots,p_k)=\frac{(p_1+p_2+\cdots+p_k)!}{p_1!p_2!\cdots p_k!}=\frac{(k-p_0)!}{p_1!p_2!\cdots p_k!}
\end{equation}
as
\begin{equation}
\mathscr{C}(p_1,p_2,\cdots,p_k)=\frac{(-1)^{k-p_0}}{\prod_{i=1}^k(i!)^{p_i}}\mathscr{M}(p_1,p_2,\cdots,p_k).
\end{equation}
With $f(x)=\sin(\pi x)$, we get
\begin{equation}
f^{k}(x)=\pi^k\sin\left(\pi x+k\frac{\pi}{2}\right).
\end{equation}
Thus, we obtain
\begin{equation}
\frac{\partial^k}{\partial x^k}\left[\frac{1}{\sin(\pi x)}\right]=\pi^k\frac{k!}{\left[\sin(\pi x)\right]^{k+1}}\sum_{\substack{p_0+p_1+p_2+\cdots +p_k=k\\p_1+2p_2+\cdots+kp_k=k}}\mathscr{C}(p_1,p_2,\cdots,p_k)\prod_{i=0}^k\left[\sin\left(\pi x+i\frac{\pi}{2}\right)\right]^{p_i}
\end{equation}
and since
\begin{equation}
\frac{\partial^k}{\partial x^k}\left[\sum_{n=-\infty}^{\infty}\frac{(-1)^n}{x+n}\right]=(-1)^kk!\sum_{n=-\infty}^{\infty}\frac{(-1)^n}{(x+n)^{k+1}}
\end{equation}
we get the final result
\begin{empheq}[box=\fbox]{align}
\pi^{k+1}=\frac{(-1)^k}{\mathscr{B}_k(x)}\sum_{n=-\infty}^{\infty}\frac{(-1)^n}{(x+n)^{k+1}}
\end{empheq}
with
\begin{empheq}[box=\fbox]{align}
\mathscr{B}_k(x)=\frac{1}{\left[\sin(\pi x)\right]^{k+1}}\sum_{\substack{p_0+p_1+p_2+\cdots +p_k=k\\p_1+2p_2+\cdots+kp_k=k}}\mathscr{C}(p_1,p_2,\cdots,p_k)\prod_{i=0}^k\left[\sin\left(\pi x+i\frac{\pi}{2}\right)\right]^{p_i}.
\end{empheq}

\section{Particular cases}

Let us first consider the case $k=0$. All the $p_i$ are equal to zero,
\begin{equation}
\mathscr{B}_0(x)=\frac{1}{\sin(\pi x)}
\end{equation}
and we recover
\begin{equation}
\pi=\sin(\pi x)\sum_{n=-\infty}^{\infty}\frac{(-1)^n}{(x+n)}.
\end{equation}
In the case $k=1$, since $p_0+p_1=1$ and $p_1=1$, we have necessarily $p_0=0$ and
\begin{equation}
\mathscr{B}_1(x)=\frac{\cos(\pi x)}{\sin^2(\pi x)}
\end{equation}
yielding
\begin{empheq}[box=\fbox]{align}
\pi^2=\frac{\sin^2(\pi x)}{\cos(\pi x)}\sum_{n=-\infty}^{\infty}\frac{(-1)^n}{(x+n)^2}.
\end{empheq}
For $k=2$, since $p_0+p_1+p_2=2$ and $p_1+2p_2=2$, we have two possibilities:

\begin{itemize}

\item $p_0=0$, $p_1=2$, $p_2=0$, yielding the term
\begin{equation}
\frac{\cos^2(\pi x)}{\sin^3(\pi x)}
\end{equation}
in the summation in $\mathscr{B}_2$ and
\item $p_0=1$, $p_1=0$, $p_2=1$, responsible for the term
\begin{equation}
\frac{1}{2\sin(\pi x)},
\end{equation}
leading to
\end{itemize}
\begin{equation}
\mathscr{B}_2(x)=\frac{1}{\sin^3(\pi x)}-\frac{1}{2\sin(\pi x)},
\end{equation}
and finally
\begin{empheq}[box=\fbox]{align}
\pi^3=\frac{2\sin^3(\pi x)}{2-\sin^2(\pi x)}\sum_{n=-\infty}^{\infty}\frac{(-1)^n}{(x+n)^3}.
\end{empheq}

\section{Product formulas for $\pi$}

Let us start with a simple proof of the so-called Euler-Wallis formula
\begin{empheq}[box=\fbox]{align}
\frac{\sin(\pi\alpha)}{\pi\alpha}=\prod_{n=1}^{\infty}\left(1-\frac{\alpha^2}{n^2}\right),
\end{empheq}
valid $\forall~\alpha\in\mathbb{R}\setminus\mathbb{Z}^*$. Setting \cite{Tissier1991}
\begin{equation}
\mathscr{S}_n(x)=\frac{1}{x+n}+\frac{1}{x-n},
\end{equation}
one has, for $n\geq 1$,
\begin{equation}
\mathscr{S}_n(x)=\frac{2x}{x^2-n^2}
\end{equation}
and
\begin{equation}
\left|\mathscr{S}_n(x)\right|\leq\frac{2x_0}{x_0^2-n^2}
\end{equation}
if $|x|\leq x_0$ and $n>x_0$ with $x_0>0$. The series
\begin{equation}
\sum_{n=1}^{\infty}\left|\mathscr{S}_n(x)\right|^2
\end{equation}
converges normally in every segment included in $\mathbb{R}\setminus\mathbb{Z}^*$. The series
\begin{equation}
\sum_{n=1}^{\infty}\ln\left|1-\frac{x^2}{n^2}\right|
\end{equation}
is therefore derivable term by term since
\begin{equation}
\frac{d}{dx}\ln\left|1-\frac{x^2}{n^2}\right|=\mathscr{S}_n(x).
\end{equation}
Subsequently, since
\begin{equation}
\pi~\mathrm{cotan}(\pi x)-\frac{1}{x}=\sum_{n=-\infty}^{\infty}\frac{1}{x+n}-\frac{1}{x}=\sum_{n=1}^{\infty}\mathscr{S}_n(x),
\end{equation}
we get
\begin{equation}
\frac{d}{dx}\sum_{n=1}^{\infty}\ln\left|1-\frac{x^2}{n^2}\right|=\sum_{n=1}^{\infty}\mathscr{S}_n(x)=\pi~\mathrm{cotan}(\pi x)-\frac{1}{x}.
\end{equation}
Setting
\begin{equation}
g(x)=\prod_{n=1}^{\infty}\left(1-\frac{x^2}{n^2}\right)
\end{equation}
and
\begin{equation}
h(x)=\frac{\sin(\pi x)}{\pi x},
\end{equation}
we have that $g$ and $h$ are proportional in each interval of $\mathbb{R}\setminus\mathbb{Z}^*$. Since
\begin{equation}
h(0)=g(0)=1,
\end{equation}
we find that $h\equiv g$ over $]-1,1[$. Let us now introduce
\begin{equation}
g_N(x)=\prod_{n=1}^N\left(1-\frac{x^2}{n^2}\right).
\end{equation}
Then
\begin{equation}
\frac{g_N(x)}{g_{N+1}(x)}\rightarrow-\frac{(1+x)}{x}
\end{equation}
which implies
\begin{equation}
x~g(x)=-(x+1)~g(x+1)
\end{equation}
and
\begin{equation}
x~h(x)=-(x+1)~h(x+1)
\end{equation}
and thus the equality of $g$ and $h$ on $]-1,1[$ implies that they are equal everywhere, and we have therefore, for any value of $x$:
\begin{equation}\label{euler-wallis}
\frac{\sin(\pi\alpha)}{\pi\alpha}=\prod_{n=1}^{\infty}\left(1-\frac{\alpha^2}{n^2}\right).
\end{equation}
The latter formula is in fact a typical example of application of the Weierstrass factorization theorem which states that every entire function can be represented as a (possibly infinite) product involving its zeroes \cite{Conway1978,Gamelin2001}. 

Let $f$ be an entire function and let $\left\{a_n\right\}$ be the nonzero zeroes of
$f$ repeated according to multiplicity. Suppose $f$ has a zero at $z=0$ of order $m\geq 0$ (where order $0$ means $f(0)\ne 0$). Then $\exists~g$ an entire function and a sequence of integers $\left\{p_n\right\}$ such that
\begin{equation}
f(z)=z^m\exp\left[g(z)\right]~\prod_{n=1}^{\infty}E_{p_n}\left(\frac{z}{a_n}\right)
\end{equation}
where $E_n(y)=1-y$ if $n=0$ and
\begin{equation}
E_n(y)=(1-y)~\exp\left(y+\frac{y^2}{2}+\cdots+\frac{y^n}{n}\right)
\end{equation}
if $n=1,2,\cdots$. It turns out that for $\sin(\pi x)$, the sequence $p_n=1$ and the function $g(z)=\log(\pi z)$ work.

Thus Euler assumed that it must be possible to represent $\sin x$ as an infinite product of linear factors given by its roots \cite{Ciaurri2015,Meyer2022}. Using this procedure and knowing that the zeroes of $\sin x$ occur at $0$, $\pm \pi$, $\pm 2\pi$, Euler derived formula (\ref{euler-wallis}). Using the Euler-Wallis product (\ref{euler-wallis}), infinite series representations for $\pi$ (or $1/\pi$) can be obtained. For instance, we get
\begin{equation}
\frac{1}{\pi}=\frac{\sqrt{2}}{4}\prod_{n=1}^{\infty}\left(1-\frac{1}{16n^2}\right)
\end{equation}
in the case $x=1/4$ and
\begin{equation}
\frac{1}{\pi}=\frac{1}{2}\prod_{n=1}^{\infty}\left(1-\frac{1}{4n^2}\right).
\end{equation}
in the case $x=1/2$. For $x=1/5$, we have
\begin{equation}
\frac{1}{\pi}=\frac{2}{5}\frac{1}{\sqrt{3-\phi}}\prod_{n=1}^{\infty}\left(1-\frac{1}{25n^2}\right)
\end{equation}
which can be put in the form
\begin{equation}
\phi=3-\frac{4\pi^2}{25}~\prod_{n=1}^{\infty}\left(1-\frac{1}{25n^2}\right)^2
\end{equation}
or equivalently
\begin{equation}
\pi^2=\frac{25}{4}(3-\phi)~\prod_{n=1}^{\infty}\left(\frac{25n^2}{25n^2-1}\right)^2,    
\end{equation}
and in the case $x=1/10$, since $\sin(\pi/10)=1/(2\phi)$, one finds
\begin{equation}
\frac{5}{\pi\phi}=\prod_{n=1}^{\infty}\left(1-\frac{1}{100n^2}\right).
\end{equation}
Similarly, for $x=1/3$, we get
\begin{equation}
\frac{3\sqrt{3}}{2\pi}=\prod_{n=1}^{\infty}\left(1-\frac{1}{9n^2}\right)
\end{equation}
and for $x=1/6$:
\begin{equation}
\frac{3}{\pi}=\prod_{n=1}^{\infty}\left(1-\frac{1}{36n^2}\right).
\end{equation}
We recall below a few other well-known infinite products, such as the Wallis formula \cite{Wallis}
\begin{equation}
\prod_{n=1}^{\infty}\left(\frac{2n}{2n-1}\right)\left(\frac{2n}{2n+1}\right)=\frac{\pi}{2}
\end{equation}
and
\begin{equation}
\prod_{n=1}^{\infty}\left(1-\frac{1}{(2n+1)^2}\right)=\frac{\pi}{4},
\end{equation}
as well as the Vi\`ete product:
\begin{equation}
\frac{\pi}{2}=\prod_{n=2}^{\infty}\frac{1}{\cos\left(\frac{\pi}{2^n}\right)}
\end{equation}
where
\begin{equation}
\cos\left(\frac{\pi}{2^n}\right)=\frac{1}{2}\underbrace{\sqrt{2+\sqrt{2+\sqrt{2+\cdots}}}}_{(n-1)~\mathrm{times}},
\end{equation}
and the two Eulerian products
\begin{equation}
\frac{\pi^2}{6}=\prod_{p~\mathrm{prime}}\frac{p^2}{p^2-1}
\end{equation}
and
\begin{equation}
\frac{\pi}{4}=\prod_{p\equiv 3~\mathrm{mod}~ 4}\left(\frac{p}{p+1}\right)\prod_{p\equiv 1~\mathrm{mod}~ 4}\left(\frac{p}{p-1}\right)=\prod_{p~\mathrm{odd~prime}}\frac{p}{p+(-1)^{\frac{p+1}{2}}}.
\end{equation}
This is of course a non-exhaustive list, and other infinite products were derived, see for instance Ref. \cite{Sondow2005,Guillera2008,Sondow2010}, or the more complicated product \cite{stack}:
\begin{equation}
\frac{\pi}{2}=\prod_{n=1}^{\infty}\left(\frac{1}{2n}\right)^{\frac{2}{2n-1}}\left[\prod_{k=1}^n\frac{(2k)^{2k}}{(2k-1)^{2k-1}}\right]^{\frac{4}{4n^2-1}}.
\end{equation}

\section{Conclusion}

We obtained, using Fourier series expansions and multiple derivative of the function $\pi/\sin(\pi x)$, series representations for positive powers of $\pi$. Leaning on the fact that the Euler-Wallis product can be derived in a straightforward manner from the same formalism, we discussed infinite products representations of $\pi$.

\appendix

\section{Proof of a relation previously used for deriving series for $\pi^{k+2}$}

In previous works \cite{Pain2022c,Pain2022d,Pain2022e}, we used the following relation attributed to Euler and mentioned as relation 13.a p. 382 in Ref. \cite{Borwein1987}. Using Eq. (\ref{cota}), we get

\begin{equation}\label{ini}
\pi~\mathrm{cotan}(\pi x)-\pi~\mathrm{cotan}(\pi a)=\sum_{n=-\infty}^{\infty}\left(\frac{1}{x-n}-\frac{1}{a-n}\right)=\sum_{n=-\infty}^{\infty}\frac{(a-x)}{(x-n)(a-n)},
\end{equation}
which is precisely the identity from which our main results in Refs. \cite{Pain2022d,Pain2022e} were derived. Of course, relation (\ref{ini}) can also be used to derive series representations of $\pi$. For instance, setting $a=1/2$ and $x=1/4$ yields the formula:
\begin{equation}
\pi=2\sum_{n=-\infty}^{\infty}\frac{1}{(2n-1)(4n-1)}.
\end{equation}

\end{document}